\documentclass[11pt]{amsart}
\usepackage{amsfonts,epsfig,fancyhdr}
\usepackage{amsmath,amssymb,amsthm, fancyhdr,latexsym,graphicx,graphics}
\usepackage{multicol,picinpar}
\usepackage{pdfsync}
\usepackage{url,hyperref}
\usepackage{euscript}
\usepackage{wrapfig}
\usepackage{caption}
\usepackage{tikz-cd}

\def \R{\mathbb{R}}

\newtheorem{theorem}{Theorem}[section]
\newtheorem{cor}[theorem]{Corollary}
\newtheorem{defi}[theorem]{Definition}

\newtheorem{rem}[theorem]{Remark}

\newtheorem{prop}[theorem]{Proposition}
\newtheorem{thm}[theorem]{Theorem}

\usepackage{url,hyperref}
\hypersetup{}

\usepackage{euscript,color}

\definecolor{red}{rgb}{1,0,0}

\hypersetup{
    bookmarks=false,         
    unicode=false,          
    pdftoolbar=true,        
    pdfmenubar=true,        
    pdffitwindow=false,     
    pdfstartview={FitH},    
    pdftitle={My title},    
    pdfauthor={Author},     
    pdfsubject={Subject},   
    pdfcreator={Creator},   
    pdfproducer={Producer}, 
    pdfkeywords={keyword1} {key2} {key3}, 
    pdfnewwindow=true,      
    colorlinks=true,       
    linkcolor=black,          
    citecolor=green,        
    filecolor=magenta,      
    urlcolor=red         
}

\def\R{\mathbb{R}}

\def\Hol{\mathrm{Hol}}

\def \oh{\hat{\omega}}
\def\Hol{\mathrm{Hol}}
\def\loc{\mathrm{loc}}

\begin{document}

\title[Cones and Cartan Geometry]
{Cones and Cartan Geometry}

\author[A. J. Di Scala]{Antonio J. Di Scala}
\author[C. E. Olmos]{Carlos E. Olmos}
\author[F. Vittone]{Francisco Vittone}


\date{\today}
\begin{abstract}
We show that the extended principal bundle of a Cartan geometry of type $(A(m,\mathbb{R}),GL(m,\mathbb{R}))$, endowed with its extended connection $\hat\omega$, is isomorphic to the principal
$A(m,\mathbb{R})$-bundle of affine frames endowed with the affine connection as defined in classical Kobayashi-Nomizu volume I.

Then we classify the local holonomy groups of the Cartan geometry canonically associated to a Riemannian manifold.
It follows that if the holonomy group of the Cartan geometry canonically associated to a Riemannian manifold
is compact then the Riemannian manifold is locally a product of cones.
\end{abstract}

\maketitle

\section{Introduction}
The main goal of this paper is to put in evidence a relation between the holonomy group of a Cartan geometry and Riemannian cones.

Any Cartan geometry $(\mathcal{P} \to M , \omega)$ has an extended principal bundle $\tilde{\mathcal{P}} \to M$ with a principal connection $\hat\omega$, and  hence a holonomy group $\Hol_x(\hat\omega)$ acting on the fiber $\tilde{\mathcal{P}}_x$ of $\tilde{\mathcal{P}}$ at $x \in M$. Our paper was motivated by the last comments following Robert Bryant's answer to MathOverflow question \cite{B16}. Namely, the question about the compactness and classification of the possible holonomy groups of the extended Cartan connection. As he writes at the end of his  answer :

{\it Probably, I should also mention that the (ordinary) holonomy of the connection $\overline{\theta}$ is what Cartan called the holonomy of (what we now call) a Cartan connection.}

 Any Riemannian manifold has a canonical Cartan geometry  associated to it as explained in \cite[p. 72, Example 1.5.1, (iii)]{CS09}
or  \cite{Sha00}.

 That the holonomy group of the extended Cartan  connection of a Riemannian manifold is not compact in general can be understood by rolling a two dimensional sphere $(S^2,g)$ with a round Riemannian metric $g$ along a maximal circle. Indeed, the parallel transport $\tau$, w.r.t. the extended connection $(\tilde{\mathcal{P}}, \hat\omega)$  along a maximal circle at $x \in S^2$, acts on the tangent space $T_x S^2$ as a translation $v \to v + t$, where $t \in T_x S^2$ has the length of the maximal circle and the direction  on which it is rolled. So $\tau$ is the generator of a non compact subgroup of the holonomy group of the Cartan extended connection of the sphere $S^2$. Less obvious is that the local holonomy group of the  extended Cartan connection of the sphere $S^2$ at a point is not compact. That is to say, even by rolling the sphere along short loops one can get any translation.

 Our second goal is to show that the extended principal bundle and the extended connection $(\tilde{\mathcal{P}}, \hat\omega)$ of the canonical Cartan geometry associated to the Riemannian manifold $(M,g)$
is essentially the affine connection defined in Kobayashi-Nomizu's book \cite[p. 127]{KN1} on the bundle of affine frames $A(M)$. Namely,

\begin{thm}\label{AffExt} Let $M$ be a $m$-dimensional manifold.
There is a 1-1 correspondence between Cartan geometries $(p: \mathcal{P} \to M, \omega)$ of type $(A(m,\mathbb{R}),GL(m,\mathbb{R}))$ on $M$ and affine connections as in Kobayashi-Nomizu's book \cite[p. 127]{KN1} on the bundle of affine frames $A(M)$. Moreover, the $A(m,\mathbb{R})$-principal bundle $A(M)$ and the extended principal bundle $\tilde{\mathcal{P}}$ are isomorphic.
The isomorphism pullbacks the extended Cartan connection $\hat\omega$ to the affine connection form $\tilde\omega_g$.
\end{thm}

 We point out that in the above mentioned Bryant's answer the extended bundle $B \times_H G \to M$ where the principal connection $\overline{\theta}$ is defined is not formally the same extended bundle $\tilde{\mathcal{P}}$ as in \v{C}ap-Slov\'ak's book. So we include in subsection \ref{sec31}, Proposition \ref{BB}, a proof that $(B \times_H G,\overline{\theta})$ and $(\tilde{\mathcal{P}},\hat\omega)$ are isomorphic as $G$-principal bundles with a connection.

Instead of using Bryant's notation we are going to use $\theta$ for the canonical form as in \cite[p. 118]{KN1} also called {\it solder form} in \cite[p. 42, 1.3.5]{CS09} (see equation (\ref{teta})). Moreover, we take the notation and definitions of Cartan geometry from
\v{C}ap-Slov\'ak's book \cite{CS09}. This is so because we think that our paper could be useful for people working  either in Parabolic Geometry, and hence familiarized with \v{C}ap-Slov\'ak's notation,  or for differential geometers that are used to Kobayashi-Nomizu's books.\\

 Our main result is the following:

\begin{thm}\label{main}
Let $(M^n,g)$ be a Riemannian manifold and let $\Hol^{\loc}_x\left(\hat\omega_{\mathcal{E}} \right)$ be its local Cartan holonomy at $x \in M$ (see Definition \ref{CartanHolR}). Let $\Hol^{\loc}_x\left(g\right)$ be the local holonomy group of the Levi-Civita connection. Let $F \times \prod_{i=1}^k M_i$ be the local De Rham decomposition of $M$ around $x \in M$, where $F$ is the local flat factor and $( M_i^{n_i},g_i)$ are irreducible.
Then we have

\begin{itemize}

\item[(i)] the affine space $A_x(M)$ splits as \[ A_x(M) = A_{x_0}(F) \times \prod_{i=1}^k A_{x_i}(M_i) \]
and  \[ \Hol^{\loc}_x\left(\oh_{\mathcal{E}}\right) =  \prod_{i=1}^k \Hol^{\loc}_{x_i}\left(\oh_i\right) \, , \]
where  $\oh_i$ is the principal connection of the Cartan geometry of $( M_i^{n_i},g_i)$, and any $h = (h_1,\cdots,h_k) \in \Hol^{\loc}_x\left(\oh_{\mathcal{E}} \right)$ acts on $v = (v_0,v_1,\cdots,v_k)$ as
\[ h \cdot v = (v_0 , h_1 \cdot v_1 , \cdots, h_k \cdot v_k) \, \, ,\]
where the map $v_i \to h_i \cdot v_i $ is an affine transformation of $A_{x_i}(M_i)$.\\

\item[(ii)] for each $i=1,\cdots,k$, the Cartan holonomy $\Hol^{\loc}_{x_i}\left(\oh_i\right)$ is isomorphic either to the Levi-Civita holonomy group $\Hol^{\loc}_{x_i}\left(g_i\right)$ or to a semidirect product $\Hol^{\loc}_{x_i}\left(g_i\right)\ltimes\R^{n_i}$.\\

\item[(iii)] $\Hol_x^{\loc}(\oh_{\mathcal{E}})$ is compact if and only if $(M^n,g)$, around $x \in M$, is a product of cones.\\

\item[(iv)] if $g$ is an Einstein (non-Ricci flat) metric, then each $\Hol^{\loc}_{x_i}\left(\oh_i\right)$ is isomorphic to $\Hol^{\loc}_{x_i}\left(g_i\right)\ltimes\R^{n_i}$, and hence $\Hol^{\loc}_x\left(\hat\omega_{\mathcal{E}}\right)$ is isomorphic to a semidirect product $\Hol^{\loc}_{x}\left(g\right)\ltimes\R^{n}$.\\

\end{itemize}

Moreover, all groups in Berger's list of holonomies of irreducible Riemannian manifolds are realizable as $\Hol(\oh_{\mathcal{E}})$ of a certain cone, with the exception of the two groups $\mathrm{SO}(2)$ and $\mathrm{Sp}(n)\cdot\mathrm{Sp}(1)$.
\end{thm}

We also obtain the following global result:

\begin{thm}\label{global} If $(M,g)$ is a complete non flat Riemannian manifold then the (restricted) Cartan holonomy group $\Hol^{0}_x\left(\hat\omega_{\mathcal{E}}\right)$ is non compact.
\end{thm}

\vspace{1cm}

\noindent
{\bf Acknowledgments:}
A.J. Di Scala wants to thanks Guilherme Machado de Freitas with whom he discussed
the main ideas of this paper in 2016. A. J. Di Scala is a member of:  GNSAGA of INdAM,
CrypTO, the group of Cryptography and Number Theory of Politecnico di Torino and of DISMA Dipartimento di Eccellenza MIUR 2018-2022.
C. Olmos is supported by FaMAF and CIEM-CONICET. F. Vittone is supported by FCEIA-UNR and CONICET.\\

\section{Basic facts of connection theory}

In the following subsections we collect some basic facts about connections on principal bundles.
We will follow the notation and definitions from \cite{KN1}. We will not develop this basic content, but we indicate precise references for the interested reader.

\subsection{Linear connections and local holonomy groups}\label{PriCon}

Let $\Gamma$ be a connection on a $G$-principal bundle $p:P \to M$. Let us indicate with $\omega_{\Gamma}$ the connection form of the connection $\Gamma$.  For each $x\in M$, the holonomy group $\Phi(x)$ and the restricted holonomy group $\Phi^0(x)$ are defined in \cite[p. 71]{KN1} as the set of isomorphisms of the fiber $p^{-1}(x)$  which come from parallel transport along loops (null-homotopic loops for $\Phi^0(x)$) in $M$ based at $x$. By fixing $p_0 \in p^{-1}(x)$ the holonomy group is realized as a subgroup of the structure group $G$ \cite[p. 72]{KN1}. Notice that by changing $p_0$ in the fiber one obtains conjugated subgroups of $G$.

The local holonomy  group $\Phi^*(x)$ is defined as the intersection of the restricted holonomy groups $\Phi^0(x,U)$, where $U$ runs through all connected neighborhoods of the point $x$, see \cite[p. 94]{KN1}.

To keep track of the connection $\Gamma$ we will indicate the holonomy group $\Phi(x)$ as $\Hol_x(\omega_{\Gamma})$, where $\omega_{\Gamma}$ is the connection form of $\Gamma$ as mentioned above. In the same way, we will indicate by $\Hol^{0}_x(\omega_{\Gamma})$ and $\Hol^{\loc}_x(\omega_{\Gamma})$ the restricted and the local holonomy groups of $\Gamma$, instead of $\Phi^0(x)$ and $\Phi^*(x)$ in Kobayashi-Nomizu notation. \\

Let now $\pi : L(M) \to M$ be the bundle of linear frames of the manifold $M$.
The canonical form $\theta$ of $L(M)$ defined by
\begin{equation}\label{teta}
\theta(X) = u^{-1}(\pi(X)) \, \, \text{for} \, X \in T_u L(M) ,
\end{equation}
where the linear frame $u$ at $\pi(u) \in M$ is considered as a linear mapping of $\R^n$ onto $T_{\pi(u)}M$  (cf. \cite[Example 5.1, p. 55]{KN1}).

Let $(M,g)$ be a Riemannian manifold of dimension $n$ and let $\nabla$ be its Levi-Civita connection.
The covariant derivative $\nabla$ is associated to a linear connection $\Gamma$ on the bundle of linear frames $L(M)$ as explained in \cite[p. 113, Chapter III]{KN1}. Let $\omega_{\Gamma}$ be the connection form of $\Gamma$.

The holonomy group $\Hol_x(\omega_{\Gamma})$ (resp.  $\Hol^{0}_x(\omega_{\Gamma})$, $\Hol^{\loc}_x(\omega_{\Gamma})$) of the Levi-Civita connection $\nabla$ can be realized as a group of orthogonal transformations $\Hol_x(g)$ (resp.  $\Hol^0_x(g)$, $\Hol^{\loc}_x(g)$) of the tangent space $T_xM$ (cf. \cite[$\S$4, Chapter II]{KN1}).\\

\subsection{Local de Rham Decomposition Theorem.} \label{sec22}

 Recall that for each $x\in M$ there is a connected neighborhood $W$ of $x$ such that $\Hol^{loc}_x(g)=\Hol^0_x(g_{U})$, for each $U\subset W$, where $\Hol^0_x(g_U)$ is the restricted holonomy group, at $x$, of the Riemannian (connected) open submanifold $U$ of $M$ (cf. \cite[Proposition 10.1, p. 95]{KN1}). In particular, one can choose $U$ to be simply connected.  Then one can apply the local de Rham decomposition Theorem to $U$ and obtain an analogous local decomposition theorem of $M$ with respect to the local holonomy group (cf. \cite[Theorem 5.4, Chapter IV]{KN1}). More precisely:

 \begin{thm}[Local de Rham decomposition Theorem] \label{localdeR} Let $(M,g)$ be a Riemannian manifold, $x\in M$ and let $T_x M=T^{(0)}_x\oplus T^{(1)}_x\oplus\cdots \oplus T^{(k)}_x$ be a decomposition into inviarnat and irreducible subspaces of $T_x M$ w.r.t. the action of the local holonomy group $\Hol^{loc}(g)$, with $T^{(0)}_x$ the set of fixed points of this action. Let $U$ be an open neighborhood of $x$ such that $\Hol^{loc}(g)=\Hol^{0}(g_U)$ and let $T^{(i)}$ the involutive distribution of $U$ obtained by parallel displacement (in $U$) of $T^{(i)}_x$. For each $i=0,\cdots, k$, let $M_i$ be the maximal integral manifold of $T^{(i)}$ through $x$.
 Then
 \begin{enumerate}
 \item there is a neighborhood $V$ of $x$ such that $V=V_0\times V_1\cdots \times V_k$ where each $V_i$ is an open neighborhood of $x$ in $M_i$ and the Riemannian metric in $V$ is the direct product of the Riemannian metrics in the $V_i$'s.
 \item The manifold $M_0$ is locally isometric to a $n_0$-dimensional Euclidean space, with $n_0=\dim M_0$.
 \item $\Hol^{loc}(g)$ is the direct product $\Hol^{loc}(g_0)\times \cdots \times \Hol^{loc}(g_k)$, of normal subgroups, where $\Hol^{loc}(g_i)$ acts trivially on $T^{(j)}_x$ if $j\neq i$ and is irreducible on $T^{(i)}_x$ for each $i=1,\cdots,k$, and $\Hol^{loc}(g_0)$ consists of the identity only.
 \end{enumerate}
 \end{thm}

\subsection{Affine connections}\label{AffinThings}

Let $M$ be a manifold of dimension $m$. An affine frame $(p, u) $ at $x \in M$ consists of a point $p \in A_x(M)$, where $A_x(M)$ is the affine space at $x \in M$ i.e. $T_x M$ regarded as an affine space, and a linear frame $u$ of the tangent space $T_x M \in M$. According to \cite[p. 127]{KN1} a generalized affine connection of $M$ is a connection $\Gamma$ in the $A(m,\mathbb{R})$-principal bundle $A(M)$ of affine frames of $M$, where $A(m,\mathbb{R})$ is the group of affine transformations of $\mathbb{R}^m$. We will follow the notation of \cite{KN1}, considering the elements of $\mathbb{R}^n$ as column vectors. So an element of $A(m,\mathbb{R})$ can be represented as
$$\left(\begin{array}{cc}
A & b\\
0 & 1\end{array}\right)$$ with $A\in GL(m,\mathbb{R})$ and $b\in \mathbb{R}^m$.
 The right action of $\left(\begin{array}{cc}
A & b\\
0 & 1\end{array}\right)$
on the affine frame $(p,u)$ is given by
\begin{equation}\label{principalAffine} (p,u) \cdot \left(\begin{array}{cc}
A & b\\
0 & 1\end{array}\right) = (p + u(b), u \circ A)
\end{equation}

\begin{rem}
In the above description of the elements of $A(m,\mathbb{R})$, we are identifying the affine space $\mathbb{R}^m$ with the hyperplane $x_{m+1}=1$ in $\mathbb{R}^{m+1}$. Notice that in \cite[p. 43]{CS09} the authors identify $\mathbb{R}^m$ with the hyperplane $x_1=1$ and hence give a different description of $A(m,\mathbb{R})$.
\end{rem}

Any generalized affine connection determines a $\mathbb{R}^n$-valued one form $\varphi$ on $L(M)$ (see \cite[Proposition 3.1, p. 127]{KN1}).
An affine connection of $M$ is a generalized affine connection whose $1$-form $\varphi$
is the canonical form $\theta$ \cite[p. 129]{KN1}.\\

\begin{defi}\label{AffHol} Let $(M,g)$ be a Riemannian manifold. Let $\Gamma_g$ (resp. $\tilde\Gamma_g$) be the principal connection in the bundle $L(M)$ of linear frames on $M$ (resp. the bundle $A(M)$ of affine frames) induced by the Levi-Civita connection of $(M,g)$. Let $\omega_{g}$ (resp. $\tilde\omega_g$) be the connection form of $\Gamma_g$ (resp. $\tilde\Gamma_g$).  The affine holonomy group $\Hol(\tilde\omega_g)_x$ of $(M,g)$ at $x \in M$ is the holonomy group of the affine connection $\tilde\Gamma_g$ on $A(M)$. The holonomy group of the Levi-Civita connection at $x \in M$  is denoted by $\Hol(g)_x$.
\end{defi}

As explained in \cite[p. 127]{KN1} the above 1-forms $\omega_g$, $\tilde\omega_g$ and $\theta$ are related by:
\begin{equation}\label{AffCoForm} \tilde\omega_g = \begin{pmatrix} \beta^* \omega_g & \beta^* \theta \\ 0 & 0 \end{pmatrix} \,  \end{equation}
where $\beta: A(M) \to L(M)$ maps $(p, u) \to u$.

This implies that $\Hol(\tilde\omega_g)_x$ is a subgroup of the semidirect product $\Hol(g)_x \ltimes T_xM$ and that the
projection $\pi : \Hol(\tilde\omega_g)_x  \to \Hol(g)_x$ is surjective.
So we have the following proposition:

\begin{prop}\label{inAffSem} The affine holonomy group $\Hol(\tilde\omega_g)_x$ is a subgroup of the semidirect product $\Hol(g)_x \ltimes T_xM$. The projection $\pi : \Hol(\tilde\omega_g)_x  \to \Hol(g)_x$ is surjective.
\end{prop}

\section{Basic facts of Cartan geometry and its holonomy group}

In the following subsections we collect basic facts about Cartan geometry.
We follow the notation and definitions from \cite[p. 71, 1.5.1]{CS09}.

A {\it Cartan geometry} of type $(G,H)$ on a manifold $M$ is a principal fiber bundle $p:\mathcal{P} \to M$ with structure group $H$, endowed with a $\mathfrak{g}$-valued $1$-form $\omega \in \Omega^1(\mathcal{P},\mathfrak{g})$ which satisfies:
\begin{enumerate}
\item the linear map $\omega_u:T_u\mathcal{P}\to \mathfrak{g}$ is a linear isomorphism for each $u\in \mathcal{P}$;
\item $(R_h)^*\omega=Ad(h^{-1})\omega$ for each $h\in H$, where $R_h$ is the right action of $h$ on $\mathcal{P}$;
\item $\omega(A^*)=A$ for all $A\in \mathfrak{h} \subset \mathfrak{g}$.
\end{enumerate}

 The notation $(\mathcal{P} \to M, \omega)$
is also used to indicate a Cartan geometry and the $1$-form $\omega$ is dubbed Cartan connection.\\

\subsection{The Cartan holonomy group of a Cartan geometry} \label{sec31}

For any Cartan geometry of type $(G,H)$ on $M$, $(p:\mathcal{P} \to M, \omega)$, there is a natural principal connection on the extended $G$-principal bundle $\tilde{p} : \tilde{\mathcal{P}}:= \mathcal{P} \times_H G \to M$. The connection form is denoted by $\hat\omega$ in \cite[p. 1040, subsection 2.2]{CGH14} and it is formally defined in \cite[p. 83]{CS09} by using \cite[Theorem 1.5.6]{CS09}.

 Notice that $\mathcal{P} \times_H G$ is the fiber bundle associated with the $H$-principal fiber bundle $\mathcal{P} \to M$ defined in \cite[p. 55]{KN1}. That is to say, the quotient of the product $\mathcal{P} \times G$ by the right action of $H$ given by $(p,g)\cdot h = (p \cdot h, h^{-1}g)$.  Then it is straightforward to check that the $G$-right action $(p,g_1) \cdot g=(p,g_1\cdot g)$ on the product $\mathcal{P} \times G$ pass to the quotient $\mathcal{P} \times_H G$. Thus $\tilde{\mathcal{P}}$ is a $G$-principal bundle.

  If one considers the map $j:\mathcal{P}\to \tilde{\mathcal{P}}$ given by $j(p)=[p,e]$, where $[p,e]$ denotes the class of $(p,e)$ in the quotient $\mathcal{P}\times_{H}G$ and $e$ is the identity of $G$, then $\hat{\omega}$ is the connection form of the only principal connection on $\tilde{\mathcal{P}}$ such that $j^{*}\hat{\omega}=\omega$  \cite[Theorem 1.5.6, p. 81]{CS09}.

 The quotient $B \times_H G$ in Bryant's answer \cite{B16}  (where $B=\mathcal{P}$) is not formally constructed as $\mathcal{P} \times_H G$ above. Indeed, Bryant considers the action of $h \in H$ on a pair $(p,g)$ given by $(p \cdot h, g h)$ and the right $G$ action on $B \times_H G$ given by $[p,g]\ast g'=[p,(g')^{-1}g]$. Anyway it is not difficult to see that the map $\Psi:\tilde{\mathcal{P}}\to B \times_H G$ given by $\Psi([p,g]):=[p,g^{-1}]$ (where the equivalence classes correspond to the respective actions of $H$) is a $G$-principal bundle isomorphism.

Consider now Bryant's $1$-form $\overline{\theta}$ of the principal connection on $B \times_H G$ associated to $\omega$. Then the following computation shows that $j^{*}(\Psi^{*}\overline{\theta})=\omega$, and hence $\Psi^{*}\overline{\theta}=\hat{\omega}$.

 Let $\rho:B\times G \to B \times_H G $ be Bryant's projection to the quotient by the diagonal $H$-action. Notice that for all $b \in B = \mathcal{P}$ \begin{equation}\label{BBeq} (\Psi \circ j)(b) = \rho(b,e) \, .\end{equation}
Let $\theta:B\times G\to \mathfrak{g}$ Bryant's modified 'difference' 1-form given by \[ \theta_{(b,g)}(X,Y)=Ad(g)(\omega_b(X)-\gamma_g(Y)). \]
So $\rho^{*}\overline{\theta}=\theta$ i.e. $$\overline{\theta}_{[b,g]}(d{\rho}_{(b,g)}(X,Y))=\theta_{(b,g)}(X,Y) \, .$$

Then for $X\in T_b\mathcal{P}$ we have:
\begin{eqnarray*}
j^{*}(\Psi^{*}\overline{\theta})(X)&=& \overline{\theta}(d\Psi(dj(X))=\overline{\theta}(d (\Psi \circ j) (X))\\
&=&\overline{\theta}(d\rho(X,0)) \, \, \text{by the above equation (\ref{BBeq})}\\
&=&\theta_{(b,e)}(X,0)\\
&=&\omega(X).
\end{eqnarray*}

So we have the following proposition:

\begin{prop} \label{BB} Let $(B \times_H G, \overline{\theta})$ be the extended principal bundle endowed with the connection $\overline{\theta}$ as in Bryant's post \cite{B16} and let $(\tilde{\mathcal{P}},\hat{\omega})$ the extended bundle as in \v{C}ap-Slov\'ak's book. Then $(B \times_H G, \overline{\theta})$ and $(\tilde{\mathcal{P}},\hat{\omega})$ are isomorphic as $G$-principal bundles with connection. That is to say, the above diffeomorphism $\Psi: \tilde{\mathcal{P}} \to B \times_H G$ is $G$-equivariant and
\[ \Psi^* \overline{\theta} = \hat{\omega} \]
\end{prop}

\begin{defi}\label{CartanHol} The holonomy group $\Hol(\hat\omega)_x$ at $x \in M$ of the Cartan geometry $(G,H)$ is the holonomy group of connection $\Gamma$ on the extended bundle $\tilde{\mathcal{P}}$ whose connection form is $\hat\omega$ as explained in \ref{PriCon}. The local and restricted holonomy groups are denoted by $\Hol^{\loc}(\hat\omega)_x$ and $\Hol^0(\hat\omega)_x$.
\end{defi}

\subsection{The Cartan geometry associated to a Riemannian manifold}

Any Riemannian manifold $(M,g)$ has a canonical Cartan geometry of type $(\mathrm{Euc}(m),O(m))$
associated to it as explained in \cite[p. 72, Example 1.5.1]{CS09} or \cite{Sha00}.
Let us explain it in detail. First of all we have a Cartan geometry of type $(A(m,R),GL(m,R))$, with Cartan connection $\omega$, due to the Levi-Civita connection $\nabla$ of $g$. Explicitelly, on the bundle $\mathcal{P}=L(M)$ the $1$-form $\omega$ is defined as
\begin{equation}\label{CartanConnection} \omega = \begin{pmatrix} \omega_g & \theta \\ 0 & 0 \end{pmatrix}
\end{equation}
where $\omega_g$ is the connection form of the Levi-Civita connection of $g$ and $\theta$ the canonical or solder form.\\

By using the metric $g$ we have $\mathcal{E} \subset \mathcal{P}$ the subset of orthonormal affine frames.
Then we restrict $p:\mathcal{P} \to M$ to $p:\mathcal{E} \to M$. The restriction $\omega_{\mathcal{E}}$ of the Cartan connection $\omega$ to $\mathcal{E}$ belongs to $\Omega^1({\mathcal{E}},\mathfrak{euc}(m))$ and satisfies
\begin{enumerate}
\item the linear map ${\omega_{\mathcal{E}}}_u:T_u\mathcal{E}\to \mathfrak{euc}(m)$ is a linear isomorphism for each $u\in \mathcal{E}$;
\item $(R_h)^*\omega_{\mathcal{E}}=Ad(h^{-1})\omega_{\mathcal{E}}$ for each $h\in O(m)$;
\item $\omega_{\mathcal{E}}(A^*)=A$ for all $A\in \mathfrak{so}(m) \subset \mathfrak{euc}(m)$.
\end{enumerate}

\subsection{The Cartan holonomy group of a Riemannian manifold $(M,g)$}

Let $(M,g)$ be a Riemannian manifold and let $(p: \mathcal{E} \to M,\omega)$ be its canonical Cartan geometry of type $(\mathrm{Euc}(m),O(m))$. That is to say, the bundle $\mathcal{E} \to M$ of orthonormal frames is a subbundle of the bundle of frames $\mathcal{P} \to M$. Hence, by construction, the principal extended bundle $ \tilde{\mathcal{E}} \to M$ is a subbundle of the principal extended bundle $\tilde{\mathcal{P}} \to M$. Moreover, again by construction, the inclusion $i: \tilde{\mathcal{E}} \to \tilde{\mathcal{P}}$ preserves the principal connections $\hat\omega$ and $\hat\omega_{\mathcal{E}}$:
\[ i^*(\hat\omega) = \hat\omega_{\mathcal{E}} \, .\]
As a consequence the horizontal distributions on $ \tilde{\mathcal{E}} \to M$ and $\tilde{\mathcal{P}} \to M$ are the same.
Indeed, they are the kernels of $\hat\omega$ and $\hat\omega_{\mathcal{E}}$ and have the same rank $m=\dim(M)$.

\begin{defi}\label{CartanHolR} The Cartan holonomy group of a Riemannian manifold
$(M,g)$ is the holonomy group $\Hol(\hat\omega_{\mathcal{E}})_x$ at $x \in M$ of the Cartan geometry
$(\mathrm{Euc}(m),O(m))$. The local and restricted holonomies
are denoted by $\Hol^{\loc}(\hat\omega_{\mathcal{E}})_x$ and $\Hol^0(\hat\omega_{\mathcal{E}})_x$.
\end{defi}

From the previous paragraph and definition we get the following proposition.

\begin{prop} \label{AffEuc} The fibers $\tilde{\mathcal{E}}_x$ of the extended bundle $\tilde{\mathcal{E}} \to M$ are invariant by the Cartan holonomy
group $\Hol(\hat\omega)_x$ of the Cartan geometry of type $(A(m,\mathbb{R}),GL(m,\mathbb{R}))$. The restriction of
$\Hol(\hat\omega)_x$ to a fiber $\tilde{\mathcal{E}}_x$ is the holonomy group $\Hol(\hat\omega_{\mathcal{E}})_x$ at $x \in M$ of the Cartan geometry $(\mathrm{Euc}(m),O(m))$.
\end{prop}

\section{ Some facts about the affine holonomy of a Riemannian manifold}

\subsection{The affine holonomy group of a product of Riemannian manifolds}
Let $(M_1,g_1)$ and $(M_2,g_2)$ be two Riemannian manifolds. Let $\tilde\omega_{g_1},\tilde\omega_{g_2}$ be their affine connection forms as in Definition \ref{AffHol}.
The goal of this section is to prove the following theorem.
\begin{thm}\label{products}  Let $x=(x_1,x_2) \in M_1 \times M_2$. The affine space $A_{x}(M_1 \times M_2)$ splits as
\[ A_{x}(M_1 \times M_2) = A_{x_1}(M_1) \times A_{x_2}(M_2) \]
and
\[ \Hol_x(\tilde\omega_g) \approx \Hol_{x_1}(\tilde\omega_{g_1}) \times \Hol_{x_2}(\tilde\omega_{g_2}) \]
where $\tilde\omega_g$ is the affine connection form of the Riemannian product: \[(M_1 \times M_2, g := g_1 \times g_2) \, .\]

Any $h=(h_1,h_2) \in  \Hol_x(\tilde\omega) $ acts on $v=(v_1,v_2)$ as \[ h \cdot v = (h_1 \cdot v_1, h_2 \cdot v_2) \]
where the maps $v_i \to h_i \cdot v_i$ , $i=1,2$, are affine transformations of $A_{x_i}(M_i)$.
\end{thm}
\it Proof. \rm  The product of bundles $L(M_1) \times L(M_2) \to M_1 \times M_2$ is a subbundle of $L(M_1 \times M_2) \to M_1 \times M_2$. Since the metric $g=g_1 \times g_2$ is a product, the Levi-Civita connection form $\omega_g$ restricted to the subbundle $L(M_1) \times L(M_2)$ splits as \[ \omega_g |_{L(M_1) \times L(M_2)} = \begin{pmatrix} \omega_{g_1} & 0 \\ 0 & \omega_{g_2} \end{pmatrix} \, . \]

In a similar way the product $A(M_1) \times A(M_2) \to M_1 \times M_2$ is a subbundle of $A(M_1 \times M_2) \to M_1 \times M_2$ and the restriction of $\tilde\omega_{g}$ to $A(M_1) \times A(M_2)$ is
  \[\tilde\omega_{g}|_{A(M_1) \times A(M_2)} =
  \begin{pmatrix}
  \beta_1^* \omega_{g_1} & 0                     & \beta_1^*\theta_1 \\
   0                     & \beta_2^*\omega_{g_2} & \beta_2^*\theta_2 \\
   0                     & 0                     &   0
   \end{pmatrix} \, ,\]
  where $\theta_i$ (resp. $\beta_i$) are the canonical forms of $M_i$ (resp. the projections $A(M_i) \to L(M_i)$) , $i=1,2$.

The above decomposition shows that the horizontal distribution $ker(\tilde\omega_{g})$ is tangent to the subbundle $A(M_1) \times A(M_2)$. That is to say, for all $(a_1,a_2) \in A(M_1) \times A(M_2)$: \[ ker(\tilde\omega_{g})_{(a_1,a_2)} \subset T_{(a_1,a_2)} A(M_1) \times A(M_2) \, .\]

So any horizontal lift of $\gamma(t)$ starting at a point $(a_1,a_2) \in A(M_1) \times A(M_2)$
is the product of the horizontal lifts of $\gamma_1, \gamma_2$ starting at $a_1, a_2$.
It follows that \[ \Hol_x(\tilde\omega_g) \approx \Hol_{x_1}(\tilde\omega_{g_1}) \times \Hol_{x_2}(\tilde\omega_{g_2}) \]
since by fixing a point $(a_1,a_2) \in A(M_1) \times A(M_2) \subset A(M_1 \times M_2)$ both groups are realized as the same subgroup of the product $A(n_1,\mathbb{R}) \times A(n_2,\mathbb{R})$.

From this it is clear that any $h=(h_1,h_2) \in  \Hol_x(\tilde\omega) $ acts on \[ v=(v_1,v_2) \in A_{x_1}M_1 \times A_{x_2}M_2 = A_{(x_1,x_2)}(M_1 \times M_2)\] as \[ h \cdot v = (h_1 \cdot v_1, h_2 \cdot v_2) \]
where the maps $v_i \to h_i \cdot v_i$ , $i=1,2$, are affine transformations of $A_{x_i}(M_i)$. $\Box$

Taking into account that the holonomy group of any principal connection on a principal bundle over $\mathbb{R}$ is trivial we get the following corollary:

\begin{cor} \label{flatfactor} Let $(\mathbb{R}^n,g_{can})$ where $g_{can} = \sum_{i=1}^n (d x_i)^2$ is the canonical metric of $\mathbb{R}^n$ and fix $x \in \mathbb{R}^n$. Then $\Hol_x(\tilde\omega_{g_{can}})$ is trivial.
\end{cor}

\subsection{The affine holonomy group of an irreducible Riemannian manifold}
Let $(M,g)$ be a Riemannian manifold.
As stated in Proposition \ref{inAffSem} $\Hol(\tilde\omega_g)_x$ is a subgroup of the semidirect product $\Hol(g)_x \ltimes T_xM$ and the projection $\pi$ to the linear part is surjective.

\begin{prop} \label{Irredu} Assume that the Riemannian manifold is irreducible i.e. $\Hol(g)_x$ acts irreducibly on $T_xM$.
Then either $\pi$ is injective or $ker(\pi) \cong T_x M$. Thus we have that either $\Hol(\tilde\omega_g)_x$ is isomorphic to $\Hol(g)_x$ or \[\Hol(\tilde\omega_g)_x \cong \Hol(g)_x \ltimes T_xM \, . \]
\end{prop}
\it Proof. \rm Assume that $\pi$ is not injective and let $\tau \neq 1 \in \ker(\pi)$. Notice that we can regard $\ker(\pi)$
as a Lie subgroup of $T_xM$. Namely, as a vector subspace $\mathbb{V}:=\ker(\pi) \subset T_xM$.\\

Let $h \in \Hol(\tilde\omega_g)_x$ be any affine parallel transport. Then
\[ h = \begin{pmatrix} \sigma & u \\ 0 & 1 \end{pmatrix}  \hspace{1cm} \tau = \begin{pmatrix} 1 & v \\ 0 & 1 \end{pmatrix}\]
where $u,v \in T_xM$, $v \neq 0$ and $\sigma \in \Hol(g)_x$. Since $\ker(\pi)$ is a normal subgroup we have that the composition $h \cdot \tau \cdot h^{-1}$ belongs to $\ker(\pi)$. Explicitly:
\[ h \cdot \tau \cdot h^{-1} =
\begin{pmatrix} \sigma & u \\ 0 & 1 \end{pmatrix} \cdot
\begin{pmatrix} 1 & v \\ 0 & 1 \end{pmatrix} \cdot
\begin{pmatrix} \sigma^{-1} &  - \sigma^{-1}u \\ 0 & 1 \end{pmatrix} =
\begin{pmatrix} 1 & \sigma v \\ 0 & 1 \end{pmatrix} \, . \]
The above equality shows that $\mathbb{V}$ is $\Hol(g)_x$-invariant. Since $v \neq 0 \in \mathbb{V}$ and  $\Hol(g)_x$ acts irreducibly we get that $\ker(\pi) = T_xM$. \qed

\section{Riemannian cones and compact affine holonomies}

Let $(M,g)$ be a Riemannian manifold. We say that $(M,g)$ is a cone at $x \in M$ if there are coordinates  $(r,x^1,\cdots, x^l)$, $r > 0 $ around $x$ such that:\\

1) the coordinate $r$ of $x$ is $r_0 > 0$ ;\\

2) the metric $g =:ds^2$ around $x$ is given by
$$ds^2=dr^2+r^2\sum f_{ij}dx^idx^j$$ where $f_{ij}$ depend only on $x^1,\cdots,x^l$.\\

\begin{prop}\label{coneStr} The Riemannian manifold $(M,g)$ is a cone at $x \in M$ if and only if there is a vector field $V$ locally defined around $x$, $V(x) \neq 0$, such that \[ \nabla_X V + X = 0 \]
for all vector fields $X$ locally defined around $x$, where $\nabla$ is the Levi-Civita covariant derivative of $g$.
\end{prop}

\it Proof .  \rm If $(M,g)$ is a cone at $x$ then a straightforward computation shows that $V := -r \frac{\partial}{\partial r} \neq 0$ at $x$ satisfies $\nabla_X V + X = 0$ for all vector fields $X$ locally defined around $x$.

For the converse let $U$ be neighborhood of $x$ where $V$ is defined and consider the function $f:U \to \mathbb{R}$, given by $f(q)=\|V(q)\|^2$. Consider the level sets of $f$ $$M^s:=\{q\in U\ :\ \|V(q)\|^2=s\}.$$
Let $F_t$ be the flow of $V$.

The following facts are straightforward computations:

\begin{itemize}

\item[(a)] the gradient of $f$ w.r.t. $g$ is $-2V$ i.e. $\mathrm{grad}(f) = -2V$.

\item[(b)] $F_t(M^s) = M^{s \cdot e^{-2t}}$

\item[(c)] If $x_1,\cdots x_n$ are local coordinates of the level set $M^s$, $s \neq 0$, then $(t,x_1,\cdots,x_n) \to F_t(x_1,\cdots,x_n)$ are local coordinates on an open subset of $M$ such that \[ g = ds^2 = s \cdot e^{-2t} dt^2+ e^{-2t} \sum f_{ij}dx^idx^j \]
    where the functions $f_{ij}$ do not depends on $t$.

\end{itemize}

Set $s=f(x)$. Then using the above chart $(t,x_1,\cdots, x_n)$ and changing $t$ by $r = -\sqrt{s} \cdot e^{-t}$ we have local coordinates around $p$ such that
\begin{equation}\label{eqCone} g = ds^2 = dr^2+ r^2 \sum s \cdot f_{ij}dx^idx^j. \end{equation}
So $(M,g)$ is a cone at $x \in M$. \qed\\

\begin{rem}\label{homothetic} Observe that the above fact {\rm (c)} on the flow $F_t$ of the vector field $V$ can be written as
\[ F_t^* g = e^{-2t} g \, ,\]
which means that the local transformation $F_t$ is homothetic but not a local isometry, cf. \cite[p. 242]{KN1}.
\end{rem}

\begin{rem}\label{kerCur} Observe that a vector field $V$ that satisfies $\nabla_X V + X = 0$ is tangent to the nullity distribution of the curvature tensor $R_{X,Y}$ of $\nabla$. That is to say, $R_{X,Y}V = 0$ for all $X,Y$.
In particular, the Ricci tensor is zero in the $V$-direction.
\end{rem}

\begin{cor}\label{fixedpoint} If the Riemannian manifold $(M,g)$ is a cone at $x \in M$ then the local holonomy group $\Hol_x^{\loc}(\tilde\omega_g)$ leaves a point of the affine tangent space $A_x(M)$ fixed.
\end{cor}
\it Proof. \rm As observed above the vector field $V := -r \frac{\partial}{\partial r}$ defined around $x$ satisfies $\nabla_X V + X = 0$ for all vector fields $X$ locally defined around $x \in M$. So by \cite[p. 195, Lemma 1 and Lemma 2]{KN1} the affine point given by $V(x) \in A_x(M)$ is a fixed point of $\Hol_x^{\loc}(\tilde\omega_g)$.\qed\\

\section{Proof of Theorem \ref{AffExt}}

The following proposition is proved inside \cite[p. 42, subsection 1.3.5]{CS09}:

\begin{prop} \label{propPLM} Let $M$ be a $m$-dimensional manifold.
There is a 1-1 correspondence between Cartan geometries $(p: \mathcal{P} \to M, \omega)$ of type
\linebreak $(A(m,\mathbb{R}),GL(m,\mathbb{R}))$ on $M$ and affine connections as in Kobayashi-Nomizu book \cite[p. 127]{KN1} on the bundle of affine frames $A(M)$. In particular, the principal bundle $p: \mathcal{P} \to M$ is $GL(m,\mathbb{R})$-isomorphic to the bundle $L(M) \to M$ of linear frames of $M$.
\end{prop}

So to prove Theorem \ref{AffExt} we have to show that the extended principal bundle $\tilde{\mathcal{P}}$
is isomorphic to the principal bundle of affine frames $A(M)$ and that the isomorphism pullbacks the extended Cartan connection $\hat\omega$ to the affine connection $\tilde\omega_g$.\\

By Proposition \ref{propPLM} we can assume that $\mathcal{P} = L(M)$.\\

We start observing that the product $L(M) \times \mathbb{R}^m$ is a $A(m,\mathbb{R})$-principal bundle over $M$ where an affine transformation $\begin{pmatrix} a & \xi \\ 0 & 1 \end{pmatrix}$ of $\mathbb{R}^m$ acts on the right as
\begin{equation}\label{prodotto} (u,v) \cdot \begin{pmatrix} a & \xi \\ 0 & 1 \end{pmatrix} = (u \circ a , a^{-1}\cdot v - a^{-1}\cdot \xi) \, \end{equation}
and $u: \mathbb{R}^m \to T_xM $ is a frame regarded as a linear map.

We are going to show that both $A(m,\mathbb{R})$-principal bundles $A(M)$ and $\tilde{\mathcal{P}}$ are isomorphic to the above product.\\

\noindent
{\bf Claim I}.  $\tilde{\mathcal{P}}$ is isomorphic, as $A(m,\mathbb{R})$-principal bundle over $M$, to the product $L(M) \times \mathbb{R}^m$.

\smallskip

The extended bundle $\tilde{\mathcal{P}}$ is $L(M) \times_{GL(m,\mathbb{R})} A(m,\mathbb{R})$ as defined in section \ref{sec31}. More precisely, $\sigma \in GL(m,\mathbb{R})$ acts on a pair $(u,\begin{pmatrix} a & \xi \\ 0 & 1 \end{pmatrix})$ as $(u \circ \sigma, \begin{pmatrix} \sigma^{-1} \cdot a & \sigma^{-1}\cdot \xi \\ 0 & 1 \end{pmatrix})$.
Let $q : L(M) \times A(m,\mathbb{R}) \to L(M) \times_{GL(m,\mathbb{R})} A(m,\mathbb{R}) $ be the projection to the quotient and let $f: L(M) \times A(m,\mathbb{R}) \to L(M) \times \mathbb{R}^m $ be the map defined as
\[ f((u,\begin{pmatrix} a & \xi \\ 0 & 1 \end{pmatrix})) = (u \circ a,  -a^{-1}\cdot \xi) \]
Form the above formulas, it is not difficult  to see that $f$ factors through the quotient projection $q$ as $f = \iota \circ q$ and by using  (\ref{prodotto})  that  $\iota:\tilde{\mathcal{P}} \to L(M) \times \mathbb{R}^m$ is a $A(m,\mathbb{R})$-equivariant diffeomorphism. Actually we have the following commutative diagram:

\[ \begin{tikzcd}
 L(M) \times A(m,\mathbb{R}) \arrow[rd,"f"] \arrow[d,"q"]\\
 \tilde{\mathcal{P}} \arrow{r}{\iota} \arrow[swap]{d}{\tilde{p}} & L(M) \times \mathbb{R}^m \arrow{d}{\pi} \\%
M \arrow{r}{id}& M
\end{tikzcd}
\]
where $x=\pi(u,v) \in M$ is the point where the linear frame $u$ is defined.

{\bf Claim II}. $A(M)$ is isomorphic, as $A(m,\mathbb{R})$-principal bundle over $M$, to the product $L(M) \times \mathbb{R}^m$. Indeed consider the map $h: A(M) \to L(M) \times \mathbb{R}^m$ defined as: \[h((p, u)) = (u , u^{-1}(o_x - p)) \]
where $p \in A_x(M)$ and $o_x$ is the zero vector of $T_xM$. Then by using (\ref{principalAffine}) and (\ref{prodotto}) it is straightforward to see that $h$ is a $A(m,\mathbb{R})$-equivariant diffeomorphism.
It follows that $\varphi: \tilde{\mathcal{P}} \to A(M)$ given by $\varphi := h^{-1} \circ \iota$ is a isomorphism of $A(m,\mathbb{R})$-principal bundles.\\

Finally, to prove that the isomorphism $\varphi$ pullbacks the affine connection $\tilde\omega_g$ to the extended Cartan connection $\hat\omega$  consider the inclusion $i: \mathcal{P} \to \tilde{\mathcal{P}} $ \[ i(u)=q(u,\mathrm{Id}) \in \tilde{\mathcal{P}} \, , \]
where $\mathrm{Id}$ is the identity of $A(m,\mathbb{R})$.\\
Observe that the composition $\beta \circ \varphi \circ i : \mathcal{P} \to \mathcal{P} $ is the identity map of $\mathcal{P}$, where $\beta: A(M) \to L(M)$ maps $(p, u) \to u$.\\
Then by using equations (\ref{AffCoForm}) and (\ref{CartanConnection}) we get
 \[ i^* \varphi^* \tilde\omega_g = i^* \varphi^* \begin{pmatrix} \beta^* \omega_g & \beta^* \theta \\ 0 & 0 \end{pmatrix} = \begin{pmatrix} i^* \varphi^* \beta^* \omega_g & i^* \varphi^* \beta^* \theta \\ 0 & 0 \end{pmatrix} = \omega \, \, . \]

So $\varphi^* \tilde\omega_g = \hat\omega$ as explained in \cite[Theorem 1.5.6, p. 81]{CS09}. \qed\\

We have the following corollary:

\begin{cor} \label{identify} Let $(M,g)$ be a Riemannian manifold. Then the Cartan holonomy groups of the canonical Cartan geometry and the affine holonomy groups at $x \in M$ of Definition \ref{AffHol} and Definition \ref{CartanHolR} are isomorphic. More precisely,
\[ \Hol_x(\hat\omega_{\mathcal{E}}) \cong \Hol_x(\tilde\omega_g) \]
\[ \Hol_x^0(\hat\omega_{\mathcal{E}}) \cong \Hol^0_x(\tilde\omega_g) \]
\[ \Hol_x^{\loc}(\hat\omega_{\mathcal{E}}) \cong \Hol^{\loc}_x(\tilde\omega_g) \]

As a consequence, the holonomy groups of the Cartan geometry at a point $x \in M$ act on the affine space $A_x(M)$ as affine transformations.
\end{cor}

\it Proof . \rm
The holonomy groups of the Cartan geometry of type $(A(m,\mathbb{R}),GL(m,\mathbb{R}))$ are by Theorem \ref{AffExt} isomorphic
to the holonomy groups of the affine connection. Then Proposition \ref{AffEuc} shows that the holonomy groups of the canonical Cartan geometry are the restriction of the holonomy groups of the Cartan geometry of type $(A(m,\mathbb{R}),GL(m,\mathbb{R}))$.  \qed

\section{Proof of Theorem \ref{main}}

Let $x\in M$.  By Theorem \ref{localdeR},  there is a neighborhood $U$ of $x$ such that
$(U,g_U)$ is isometric to a Riemannian product \[ (F,g_0) \times \prod_{i=1}^k (M_i,g_i) \]
 around $x \in M$ where $(F,g_0)$ is the flat factor which, i.e., it is isometric to an open subset of a Euclidean space. Accordingly,  the local holonomy group of the Levi-Civita connection splits as
\[ \Hol_x^{\loc}(g) = \{e\} \times \prod_{i=1}^k \Hol_{x_i}(g_i) \]
and each $\Hol_{x_i}(g_i)$ acts irreducible on $T_{x_i} M_i$.

\subsection{Proof of Theorem \ref{main}}

Part (i) follows from Theorem \ref{products} and Corollary \ref{identify} applied to the local De Rham decomposition at $x \in M$:
\[ (F,g_0) \times \prod_{i=1}^k (M_i,g_i) \]
Observe that the affine holonomy of the flat factor $(F,g_0)$ is trivial by Corollary \ref{flatfactor}.\\

Part (ii) follows from Proposition \ref{Irredu} and Corollary \ref{identify}.\\

Part (iii). Assume that $\Hol^{\loc}_x(\hat\omega_{\mathcal{E}})$ is compact. We claim that each $(M_i,g_i)$ is a cone at $x_i \in M_i$. Indeed, each $\Hol^{\loc}_{x_i}(\tilde\omega_{g_i})$ is compact by Corollary \ref{identify} and acts as affine transformations of the affine space $A_{x_i}M$. Then by \cite[p. 183, Lemma 2.3]{M77} $\Hol^{\loc}_{x_i}(\tilde\omega_{g_i})$ has a fixed point. Then by \cite[p. 195, Lemma 1 and Lemma 2]{KN1} the fixed point gives raise to a vector field $V$ defined around $x_i$ such that $\nabla_X V + X = 0$ for all vector fields $X$ defined around $x_i$.\\

If $V(x_i) \neq 0$ then Proposition \ref{coneStr} shows that $(M_i,g_i)$ is a cone at $x_i \in M_i$ and we are done.\\

If $V(x_i) = 0$ we get a contradiction. Indeed, assume $V(x_i) = 0$ and
let $f:U \to \mathbb{R}$ as in the proof of Proposition \ref{coneStr} where $U$ is a neighborhood of $x_i$.
Then the Hessian of $f$ w.r.t. $g_i$ at $x_i$ is $\mathrm{Hess}_f(X,Y)_p = g_i(\nabla_X (\mathrm{grad} f)_p,Y) = 2 g(X,Y)$.
So by Morse Lemma \cite[Lemma 2.2., p. 6]{M63} we can assume that there are coordinates $u_1,\cdots,u_l$, $l=\mathrm{dim}(M_i)$, on $U$ such that $f = \sum_{j=1}^{l} u_j^2$. The subsets $U_s = \{ q \in U:  f(q) < s \}$ are (topological) balls and the flow $F_t$ of $V$ preserves $U_s$ for all $t>0$ i.e. $F_t (U_s) \subset U_s$ for $t>0$. At all points of $U_s \setminus \{x_i\}$ there are coordinates as in equation (\ref{eqCone}) so the Lie derivative $\mathcal{L}_V g = -2 \cdot g$. So for $t>0$ fixed $F_t^* g = e^{-2t} g$ hence $F_t$ is a homothetic transformation of $U_s$. Then the proof of \cite[Lemma 2, p. 242]{KN1} shows that the restriction of $g$ to $U$ is a flat metric. Thus the local holonomy group of the Levi-Civita connection at $x_i \in M_i$ is trivial. This contradicts that $M_i$ is locally irreducible (and non flat) at $x_i$.\\

To finish Part (iii) assume that $(M,g)$ is a product of cones around $x$.
Then by using Theorem \ref{products} and Corollary \ref{fixedpoint} we get that $\Hol^{\loc}_x(\hat\omega_{\mathcal{E}})$ leaves a point of the affine space $A_x(M)$ fixed. Hence, $\Hol^{\loc}_x(\hat\omega_{\mathcal{E}})$ is compact by Part (ii).\\

Part (iv). If $g$ is Einstein (non Ricci flat) then there is no local flat factor $F$ at $x$, and each $g_i$ is an Einstein metric (non Ricci flat). By Proposition \ref{Irredu} each $\Hol^{\loc}_{x_i}\left(\oh_i\right)$ is either $\Hol^{\loc}_{x_i}\left(g_i\right)\ltimes\R^{n_i}$
or compact. If $\Hol^{\loc}_{x_i}\left(\oh_i\right)$ is compact then by Part (iii) $(M_i,g_i)$
is a cone at $x_i$. Observe that the Ricci tensor of a cone vanishes in the radial direction $\frac{\partial}{\partial r}$ by Remark \ref{kerCur}.
Since we assume the metric $g_i$ to be non Ricci flat we get that $\Hol^{\loc}_{x_i}\left(\oh_i\right)$ is not compact.
So $\Hol^{\loc}_{x_i}\left(\oh_i\right)$ is isomorphic to a semidirect product $\Hol^{\loc}_{x_i}\left(g_i\right)\ltimes\R^{n_i}$. Thus $\Hol^{\loc}_x\left(\hat\omega_{\mathcal{E}}\right)$ is isomorphic to a semidirect product $\Hol^{\loc}_{x}\left(g\right)\ltimes\R^{n}$.\\

The statement that all irreducible holonomies of Berger's list but $\mathrm{SO}(2)$ and $\mathrm{Sp}(n)\cdot\mathrm{Sp}(1)$ can be realized as holonomies of cones is proved in \cite[pag. 7]{ACL19}. That $\mathrm{Sp}(n)\cdot\mathrm{Sp}(1)$ cannot be realized as the Levi-Civita holonomy of a Riemannian cone follows because a manifold with such holonomy is Einstein non Ricci-flat and the Ricci tensor of a Riemannian cone is zero in the $\frac{\partial}{\partial r}$ direction by Remark \ref{kerCur}.  Observe that $\mathrm{SO}(2)$ can not be the holonomy of a cone since a 2-dimensional cone is flat hence the holonomy is trivial.

\section{Proof of Theorem \ref{global}}

The proof is by contradiction. Assume that the restricted Cartan holonomy group is compact.
Then the restricted affine holonomy group is compact by
Corollary \ref{identify}. So the restricted affine holonomy group has a fixed point by \cite[Lemma 2.3., p. 183]{M77}, and
hence $(M,g)$ is flat by \cite[Corollary 7.3., p. 197]{KN1}, which is a contradiction. \qed \\

\pagebreak
\noindent {\bf Authors' Addresses:}

\vspace{.5cm}

\noindent
\begin{tabular}{l|cl|cl}
Antonio J. Di Scala & & Carlos E. Olmos  & & Francisco Vittone\\
\footnotesize Dipartimento di Scienze Matematiche&  & \footnotesize FaMAF, CIEM-Conicet & &\footnotesize DM, ECEN, FCEIA - Conicet\\
\footnotesize Politecnico di Torino &  &\footnotesize Universidad Nac. de C\'ordoba& &\footnotesize Universidad Nac. de Rosario \\
 \footnotesize Corso Duca degli Abruzzi, 24& &\footnotesize Ciudad Universitaria  & &\footnotesize Av. Pellegrini 250\\
 \footnotesize 10129, Torino, Italy  & & \footnotesize 5000, C\'ordoba, Argentina & &\footnotesize 2000, Rosario, Argentina \\
\footnotesize{antonio.discala@polito.it} & & \footnotesize{olmos@famaf.unc.edu.ar} & & \footnotesize{vittone@fceia.unr.edu.ar}\\
\footnotesize{http://calvino.polito.it/$\sim$adiscala}& &
 & &\footnotesize{www.fceia.unr.edu.ar/$\sim$vittone}
\end{tabular}

\end{document}